\documentclass[12pt,reqno,centertags]{amsart}
\usepackage{amsmath,latexsym}
\usepackage[psamsfonts]{amssymb}
\usepackage{amssymb}
\usepackage{upref}
\usepackage{mathptm}
\DeclareSymbolFont{SY}{U}{psy}{m}{n}
\DeclareMathSymbol{\emptyset}{\mathord}{SY}{'306}

\newcommand{\bbR}{\mathbb R}

\newcommand{\bbC}{\mathbb C}

\newcommand{\sE}{\mathsf E}

\newcommand{\sP}{\mathsf P}
\newcommand{\sQ}{\mathsf Q}

\newcommand{\wOm}{\widetilde{\Omega}}
\renewcommand{\epsilon}{\varepsilon}
\newcommand{\lal}{\langle}
\newcommand{\ral}{\rangle}
\newcommand{\be}{\begin{equation}}
\newcommand{\ee}{\end{equation}}

\newcommand{\spec}{\sigma}
\newcommand{\R}{\mathbb{R}}
\newcommand{\EE}{\mathsf{E}}

\newcommand{\cB}{{\mathcal B}}

\newcommand{\cG}{{\mathcal G}}
\newcommand{\cH}{{\mathcal H}}

\newcommand{\cK}{{\mathcal K}}
\newcommand{\cL}{{\mathcal L}}

\newcommand{\cN}{{\mathcal N}}

\newcommand{\cQ}{{\mathcal Q}}

\newcommand{\cX}{{\mathcal X}}

\newcommand{\dist}{{\ensuremath{\mathrm{dist}}}}

\newcommand{\diag}{{\ensuremath{\mathrm{diag}}}}

\newcommand{\ri}{\mathrm{i}}
\DeclareMathOperator{\Ran}{\mathrm{Ran}}

\newcommand{\Dom}{\mathrm{Dom}}
\newcommand{\dom}{\mathrm{Dom}}

\newtheorem{introtheorem}{Theorem}{\bf}{\it}

\newtheorem{thm}{Theorem}[section]
\newtheorem{hypo}[thm]{Hypothesis}
\newtheorem{cor}[thm]{Corollary}

\theoremstyle{definition}
\newtheorem{defn}[thm]{Definition}
\theoremstyle{remark}
\newtheorem{rem}[thm]{Remark}

\numberwithin{equation}{section}

\begin{document}

\title[Solutions to the Riccati equation and the tan $\Theta$
theorem] {On the existence of solutions to the
operator
%%\newline
Riccati equation and the
tan\,$\boldsymbol{\Theta}$ theorem}

\author[V. Kostrykin]{V. Kostrykin}
\address{Fraunhofer-Institut f\"{u}r Lasertechnik \\
Steinbachstra{\ss}e 15 \\
Aachen, D-52074, Germany}
\email{kostrykin@ilt.fraunhofer.de, kostrykin@t-online.de}

\author[K. A. Makarov]{K. A. Makarov}
\address{Department of Mathematics \\
University of Missouri \\
Co\-lum\-bia, MO 65211, USA}
\email{makarov@math.missouri.edu}

\author[A. K. Motovilov]{A. K. Motovilov}
\address{Department of Mathematics \\
University of Missouri \\
Columbia, MO 65211, USA. \newline
\hspace*{12pt}{\it Permanent address: }
BLTP, JINR\\
141980 Dubna\\
Moscow Region\\
Russia} \email{motovilv@thsun1.jinr.ru}

\subjclass{Primary 47A55, 47A15; Secondary 47B15}
\keywords{Operator Riccati equation, graph subspaces,
block operator matrices}

\begin{abstract}
Let $A$ and $C$ be self-adjoint operators such that the
spectrum of $A$ lies in a gap of the spectrum of $C$ and
let $d>0$ be the distance between the spectra of $A$ and
$C$. We prove that under these assumptions the sharp value
of the constant $c$ in the condition $\|B\|<c d$
guaranteeing the existence of a (bounded) solution to the
operator Riccati equation $XA-CX+XBX=B^*$ is equal to
$\sqrt{2}$. We also prove an extension of the Davis-Kahan
$\tan\Theta$ theorem and provide a sharp estimate for the
norm of the solution to the Riccati equation.  If $C$ is
bounded, we prove, in addition, that the solution $X$ is a
strict contraction if $B$ satisfies the condition
$\|B\|<d$, and that this  condition is sharp.
\end{abstract}

\subjclass{Primary 47A55, 47A15; Secondary 47B15}

\keywords{invariant subspaces, operator Riccati equation,
operator matrix, operator angle, factorization theorem}

\maketitle

\thispagestyle{empty}

%%%%%%%%%%%%%%%%%%%%%%%%%%%%%%%%%%%%%%%%%%%%%%
\section{Introduction}
\label{Intro}
%%%%%%%%%%%%%%%%%%%%%%%%%%%%%%%%%%%%%%%%%%%%%%

In this paper we consider the operator Riccati equation
\begin{equation}
\label{Ric0} XA-CX+XBX=B^*
\end{equation}
associated with the self-adjoint $2\times 2$ block
operator matrix
\begin{equation}
\label{cHintro} {H}=\begin{pmatrix}
  A     &   B  \\
  B^*   &   C
\end{pmatrix}
\end{equation}
on the orthogonal sum $\cH=\cH_A \oplus \cH_C$ of
separable Hilbert spaces $\cH_A$ and $\cH_C$. Here $A$ is
a bounded self-adjoint operator on the Hilbert space
$\cH_A$, $C$ possibly unbounded self-adjoint operator on
the Hilbert space $\cH_C$, and $B$ is a bounded operator
from $\cH_C$ to $\cH_A$.

Solving the Riccati equation appears to be an adequate
tool in the study of the invariant subspaces of the
operator $H$ that are the graphs of bounded operators from
$\cH_A$ to $\cH_C$, the so-called graph subspaces. It is
well known that given a bounded solution $X:\cH_A\to
\cH_C$ to the Riccati equation \eqref{Ric0} (with $\Ran
X\subset \Dom(C)$), the graph
\begin{equation*}
\cG(X)=\left\{x\oplus X x|\ x\in\cH_A \right\}
\end{equation*}
of the operator $X$ reduces the operator $H$  (see, e.g.,
\cite[Section 5]{AMM}) and in the framework of this
approach the following two problems naturally arise. The
first problem is to study the spectrum of the part of $H$
associated with the reducing subspace $\cG(X)$
(respectively its orthogonal complement $\cG(X)^\perp$),
and the second one is to estimate the operator angle
$\Theta$ (see, e.g., \cite{KMMgeom} for discussion of this
notion) between the subspaces $\cG(X)$ and $\cH_A$
(respectively $\cG(X)^\perp$ and $\cH_C$). Both of these
problems can efficiently be solved if a bounded solution
$X$ to \eqref{Ric0} is known: the operator $H$ appears to
be similar to the diagonal block operator matrix
\begin{equation*}
\begin{pmatrix} A + B X & 0 \\ 0 & C- B^\ast X^\ast \end{pmatrix}
\end{equation*}
associated with the decomposition $\cH=\cH_A\oplus \cH_C$
(in particular, $\sigma(H)=\sigma(A + B X)\cup \sigma (C-
B^\ast X^\ast )$), and the operator angle $\Theta$ between
the subspaces $\cG(X)$ and $\cH_A$ has the representation
\begin{equation*}
\Theta = \arctan \sqrt{X^\ast X}.
\end{equation*}

If the spectra of $A$ and $C$ overlap, the Riccati
equation \eqref{Ric0} may  have no solution at all (cf.,
e.g., \cite[Example 3.2]{AMM}). At the same time the
spectra separation requirement alone does not guarantee
the existence of solutions either (see, e.g.,  \cite[Lemma
3.11]{AMM}). Under the spectra separation hypothesis
\begin{equation}\label{spsep}
\dist(\sigma(A), \sigma(C)) > 0,
\end{equation}
a natural sufficient condition for the existence of
solutions to the Riccati equation \eqref{Ric0} requires a
smallness assumption on the operator $B$ of the form
\begin{equation}\label{cdest}
\|B\|<c_{\text{best}}\,\dist\bigl(\spec(A),\spec(C)\bigr)
\end{equation}
with a  constant $c_{\text{best}}>0$ independent of the
distance between the spectra $\sigma(A)$ and $\sigma(C)$
of the operators $A$ and $C$, respectively. The best
possible constant $c_{\mathrm{best}}$ in \eqref{cdest} is
still unknown. However, $c_{\mathrm{best}}$ is known to be
in the interval $\left[{\pi}^{-1},\sqrt{2}\right]$ (see
\cite{AMM}). If both $A$ and $C$ are bounded, then
$c_{\mathrm{best}}\in\left[c_\pi,\sqrt{2}\right]$ with
$c_\pi=\frac{3\pi-\sqrt{\pi^2+32}}{\pi^2-4}=0.503288...$
(see \cite{KMMgamma}). In \cite{KMMgamma} the best
possible constant $c_{\mathrm{best}}$ has been conjectured
to be $\sqrt{3}/2$. Some earlier results in this direction
can be found in \cite{AdLT}, \cite{MeMo99},
\cite{Motovilov:SPb:91}, and \cite{Motovilov:95}.

In some particular cases  the optimal solvability
condition \eqref{cdest} can be relaxed provided that some
additional assumptions upon mutual disposition of the
spectra of $A$ and $C$ are posed.  For instance, if the
spectra of $A$ and $C$ are subordinated, e.g.,
\begin{equation}
\label{subord} \sup\spec(A)<\inf\spec(C),
\end{equation}
the Riccati equation \eqref{Ric0} is known to have a
strictly contractive solution for any bounded $B$ (see,
e.g., \cite{AL95}). To some extent abusing the terminology
one may say that in this case the best possible constant
in inequality \eqref{cdest} is infinite: No smallness
assumptions on $B$ are needed.

In the limiting case of \eqref{subord},
\begin{equation*}
\sup\spec(A)=\inf\spec(C),
\end{equation*}
the existence of contractive solutions has been
established in \cite{AdLT} under some additional
assumptions which have been dropped in \cite{KMMalpha}.
See also \cite{MenShk} where the spectra separation
condition \eqref{spsep}  has also  been somewhat relaxed
and the existence of a bounded but not necessarily
contractive solution has been established.

Our \emph{first principal result} concerns  the case where
the operator $C$ has a finite spectral gap  containing the
spectrum of $A$. Recall that  by a finite spectral gap of
a self-adjoint operator $T$ one  understands an
\emph{open} finite interval on the real axis lying in the
resolvent set of $T$ such that both of its end points
belong to the spectrum of $T$.

\begin{introtheorem}\label{thm:1}
Assume that the self-adjoint operator $C$ has a finite
spectral gap $\Delta$ containing the spectrum of the
bounded self-adjoint operator $A$.
\smallskip

(i) Suppose that
\begin{equation}
\label{korint} \|B\|<\sqrt{d|\Delta|} \quad \text{where}
\quad d=\dist(\spec(A),\spec(C)),
\end{equation}
with $|\cdot|$ denoting Lebesgue measure on $\bbR$. Then
the spectrum of the block operator matrix
$H=\begin{pmatrix}A & B\\
B^* & C\end{pmatrix}$ in the gap $\Delta$ is a (proper)
closed subset of (the open set) $\Delta$. The spectral
subspace of the operator $H$ associated with the interval
$\Delta$ is the graph of a bounded solution $X:\
\cH_A\rightarrow\cH_C$ to the Riccati equation
\eqref{Ric0}. Moreover, the operator $X$ is the unique
solution to the Riccati equation in the class of bounded
operators with the properties
\begin{equation}
\label{Uniq} \spec(A+B X)\subset\Delta\quad\text{and}\quad
\spec(C-B^* X^*)\subset\R\setminus\Delta.
\end{equation}
\smallskip

(ii) If in addition the operator $C$ is bounded and
\begin{equation}
\label{korints} \|B\|<\sqrt{d(|\Delta|-d)},
\end{equation}
then the solution $X$ is a strict contraction and
\begin{equation}
\label{XestFin2} \|X\|\leq \tan\frac{1}{2}
\arctan\left(\frac{2\|AB+BC\|}{d(|\Delta|-d)-\|B\|^2}\right)<1.
\end{equation}
\end{introtheorem}

As a corollary, under the assumption that the operator $C$
has a finite spectral gap $\Delta$ containing the spectrum
of $A$, we prove that $c=\sqrt{2}$ is best possible in
condition \eqref{cdest} ensuring the existence of a
bounded solution to the Riccati equation \eqref{Ric0} (see
Remark \ref{coptim}) while $c=1$ is best possible in
\eqref{cdest} to ensure that the solution is a contraction
(see Remark \ref{controptim}).

The proof of the part (i) of Theorem \ref{thm:1} will be
given in Section \ref{SecExist} and that of the part (ii)
in Section \ref{SecEst}.

Our \emph{second principal result} holds with no {\it a
priori}  assumption upon the mutual disposition of the
spectra of $A$ and $C$ (in particular, the spectra of $A$
and $C$ may overlap).

\begin{introtheorem}\label{thm:2}
Assume that the self-adjoint operator $C$ has a spectral
gap $\Delta$ (finite or infinite) and the self-adjoint
operator $A$ is bounded. Assume that the Riccati equation
\eqref{Ric0} has a bounded solution $X$ and hence the
graph subspace $\cG(X)$ reduces the block operator matrix
$H$. Suppose that the spectrum of the part $H|_{\cG(X)}$
of the operator $H$ associated with the reducing subspace
$\cG(X)$ is a (proper) closed subset of (the open set)
$\Delta$. Then the the following norm estimate holds:
\begin{equation}
\label{NXest} \|X\|\leq \frac{\|B\|}{\delta} \quad
\text{with} \,\,\,
\delta=\dist\bigl(\spec(H|_{\cG(X)}),\spec(C)\bigr).
\end{equation}
Equivalently,
\begin{equation}\label{TanTheta}
\|\tan\Theta\|\leq\frac{\|B\|}{\delta},
\end{equation}
where $\Theta$ is the operator angle between the subspaces
$\cH_A$ and $\cG(X)$.
\end{introtheorem}

Estimate \eqref{TanTheta} extends the Davis-Kahan
$\tan\Theta$ theorem \cite{Davis:Kahan}, a result
previously known only in the case where the spectra of $C$
and $H|_{\cG(X)}$ are subordinated, that is, the operator
$C$ is semibounded and the spectrum of the part
$H|_{\cG(X)}$ lies in the \emph{infinite} spectral gap of
$C$. This generalization extends the list of the
celebrated Davis-Kahan $\sin\Theta$ and $\sin 2\Theta$
theorems, proven in the case where the operator $C$ has a
gap of finite length \cite{Davis:Kahan}.

The proof of Theorem \ref{thm:2} will be given in Section
\ref{SecTan}.

Our main techniques are based on applications of the
Virozub-Matsaev factorization theorem for analytic
operator-valued functions \cite{ViMt} (in the spirit of
the work \cite{MenShk} (cf.\ \cite{LMMT})) and the
Daletsky-Krein factorization formula \cite{DK}. Under the
hypothesis of Theorem \ref{thm:1} we prove that
\begin{itemize}

\item for $\lambda\notin\spec(C)$  the operator-valued Herglotz function
$M(\lambda)=\lambda I-A+B(C-\lambda I)^{-1}B^*$ admits a
factorization
\begin{equation}
\label{ZAint} M(\lambda)=W(\lambda)(Z-\lambda I),
\end{equation}
with $W$ being an operator-valued function holomorphic on
the resolvent set of the operator $C$ and $Z$ is a bounded
operator with the spectrum in the spectral gap $\Delta$ of
the operator $C$,

\item the Riccati equation \eqref{Ric0} has
a bounded solution of the form
\begin{equation}
\label{Xintro} X=-\frac{1}{2\pi{\mathrm i}}\int_\Gamma
d\lambda(C-\lambda I)^{-1}B^* (Z-\lambda I)^{-1},
\end{equation}
where $\Gamma$ is an appropriate Jordan contour encircling
the spectrum of the operator~$Z$,

\item the spectral subspace of the $2\times 2$ block operator
matrix $H$ \eqref{cHintro} associated with the interval
$\Delta$ is the  graph of the operator $X$,

\item the spectrum of the operator $H$ in the interval
$\Delta$ coincides with that of the operator $Z$, that is,
$\spec(H)\cap\Delta=\spec(Z)$.
\end{itemize}

In Section \ref{SecFact} we recall the concept of
invariant graph subspaces for linear operators as well as
their relation to the Riccati equation. Theorem
\ref{RicSol} below presents a general result linking the
factorization property \eqref{ZAint} of the operator
valued function $M(\lambda)$ with the existence of a
spectral subspace for the $2\times2$ self-adjoint block
operator matrix \eqref{cHintro} admitting representation
as the graph of the operator \eqref{Xintro}. In
Section~\ref{SecExist} under hypothesis \eqref{korint} we
prove factorization formula \eqref{ZAint} and  give bounds
on the location of the spectrum of the operator $Z$
(Theorem \ref{Zexist}), and finally prove the part (i) of
Theorem \ref{thm:1}. The proof of Theorem \ref{thm:2} is
given in Section~\ref{SecTan}. In Section~\ref{SecEst}
combining the results of Theorems \ref{thm:1} (i) and
\ref{thm:2} under assumptions \eqref{korint} and
\eqref{korints} respectively we provide norm estimates on
the solution $X$  of the Riccati equation and prove
Theorem 1 (ii). We conclude the section by an example
showing that condition \eqref{korints} ensuring the strict
contractivity of the solution  $X$ is sharp.

Few words about notations used throughout the paper. Given
a Hilbert space $\cK$ by $I_{\cK}$ we denote the identity
operator on $\cK$. If it does not lead to any confusion we
will simply write $I$ instead of more pedantic notation
$I_\cK$. The set of all bounded linear operators from the
Hilbert space $\cK$ to a Hilbert space $\cL$ will be
denoted by $\cB(\cK,\cL)$. If $\cL=\cK$ the shorthand
$\cB(\cK)$ will be used for this set. Let $K$ and $L$ be
self-adjoint operators on a Hilbert space $\cK$. We say
$K< L$ (or, equivalently, $L> K$) if there is a number
$\gamma>0$ such that $L-K>\gamma I$. The notation
$\rho(T)$ will be used for the resolvent set of a closed
operator $T$.

After completing this work we learned that a result
similar to the part (ii) of Theorem \ref{thm:1} has been
recently obtained within a different approach by
A.\,V.~Se\-lin (private communication).

%%%%%%%%%%%%%%%%%%%%%%%%%%%%%%%%%%%%%%%%%%%%%%%%%%%%%%%%%%%%%%%%%%%%%
\section{Invariant graph subspaces and block diagonalization}
%%%%%%%%%%%%%%%%%%%%%%%%%%%%%%%%%%%%%%%%%%%%%%%%%%%%%%%%%%%%%%%%%%%%%
\label{SecFact}

In this section we collect some results related to the
invariant graph subspaces of a linear operator as well as
to the closely related problem of block diagonalization of
block operator matrices.

First,  recall the definition a graph subspace.

\begin{defn}
Let $\cK$ be a closed subspace of a Hilbert space $\cN$
and $X\in\cB(\cK,\cK^\perp)$. Denote by $P_\cK$ and
$P_{\cK^\perp}$ the orthogonal projections in $\cN$ onto
the subspace $\cK$ and orthogonal complement $\cK^\perp$,
respectively. The set
\begin{equation*}
\cG(X)=\{x\in\cN\,|\, P_{\cK^\perp}x=XP_{\cK}x\}
\end{equation*}
is called the graph subspace associated with the operator
$X$.
\end{defn}

%%%%%%%%%%%%%%%%%%%%%%%%%%%%%%%%%%%%%%%%%%%%%%%%%%%%%%%%%%%%%%%%%%%%%%%%%%
For notational setup we assume the following
%%%%%%%%%%%%%%%%%%%%%%%%%%%%%%%%%%%%%%%%%%%%%%%%%%%%%%%%%%%%%%

\begin{hypo}
\label{Hmatr} Let $H_0$ be a self-adjoint operator in a
Hilbert space $\cH$ and  $\cH_A\subset\cH$  a reducing
subspace of $H_0$. Assume that
 with respect to the decomposition
\begin{equation}
\label{decom} \cH=\cH_A\oplus\cH_C \quad
(\cH_C=\cH\ominus\cH_A)
\end{equation}
the operator $H_0$ reads  as the block diagonal operator
matrix
\begin{equation*}
H_0=\mathop{\mathrm{diag}}(A,C)
\end{equation*}
with $A$ being  a bounded self-adjoint operator in
$\cH_A$, $C$ a possibly unbounded self-adjoint operator in
$\cH_C$, and $\Dom(H_0)=\cH_A\oplus\Dom(C)$. Assume, in
addition, that with respect to the decomposition
\eqref{decom}  the self-adjoint operator $H$ reads  as
\begin{equation}
\label{twochannel} H=\begin{pmatrix}
  A    &   B \\
  B^*  &   C
\end{pmatrix},\quad\Dom(H)=\cH_A\oplus\Dom(C),
\end{equation}
where $B$ is a bounded operator from $\cH_C$ to $\cH_A$.
\end{hypo}

Under Hypothesis \ref{Hmatr}, a bounded operator $X$ from
$\cH_A$ to $\cH_C$ is said to be a solution to the Riccati
equation \eqref{Ric0} if $\Ran(X)\subset\Dom(C)$ and
\eqref{Ric0} holds as an operator equality.

The existence of a bounded solution to the Riccati
equation \eqref{Ric0} is equivalent to a possibility of
the block diagonalization of the operator matrix $H$ with
respect to the decomposition
$\cH=\cG(X)\oplus\cG(X)^\perp$. The precise statement is
as follows (see Lemma 5.3 and Theorem 5.5 in \cite{AMM};
also cf.~\cite{AdLT}, \cite{Daughtry}, \cite{KMMgeom},
\cite{Motovilov:95}).
%%%%%%%%%%%%%%%%%%%%%%%%%%%%%%%%%%%%%%%%%%%%%%%%%%%%%%%%%%%%%%%%%%%%
%%%%%%%%%%%%%%%%%%%%%%%%%%%%%%%%%%%%%%%%%%%%%%%%%%%%%%%%%%%%%%%%%%%%
\begin{thm}
\label{thHi2} Assume Hypothesis \ref{Hmatr}. Then a
bounded operator $X$ from $\cH_A$ to $\cH_C$ is a solution
to the Riccati equation \eqref{Ric0} iff the graph
$\cG(X)$ of $X$ reduces the $2\times2$ block operator
matrix $H$. Moreover, if $X\in\cB(\cH_A,\cH_C)$ is a
solution to \eqref{Ric0} then:

\medskip

\noindent{\rm(i)} The operator $V^{-1}HV$  with
\begin{equation*}
V=\begin{pmatrix}
I & -X^* \\
 X  &  I
\end{pmatrix}
\end{equation*}
is block diagonal with respect to decomposition
\eqref{decom}. Furthermore,
\begin{equation*}
V^{-1}HV=\begin{pmatrix}
Z     &   0      \\
0     &   \widehat{Z}
\end{pmatrix},
\end{equation*}
where $Z=A+BX$ with $\Dom(Z)=\cH_A$  and
$\widehat{Z}=C-B^*X^*$ with $\Dom(\widehat{Z})=\dom(C)$.
\medskip

\noindent{\rm(ii)} The operator
\begin{equation}\label{HA}
\Lambda =(I+X^*X)^{1/2}Z (I+X^*X)^{-1/2}
\end{equation}
and possibly unbounded operator
\begin{equation}
\label{HC} \widehat{\Lambda} = (I+XX^*)^{1/2}\widehat{Z}
(I+XX^*)^{-1/2}
\end{equation}
with $\dom(\widehat \Lambda)=(I+XX^*)^{1/2}(\dom(C))$ are
self-adjoint operators in $\cH_A$ and $\cH_C$,
respectively.
\end{thm}
%%%%%%%%%%%%%%%%%%%%%%%%%%%%%%%%%%%%%%%%%%%%%%%%%%%%%%%%%%%%%%%
Theorem \ref{thHi2} yields the following uniqueness result
as a corollary.
%%%%%%%%%%%%%%%%%%%%%%%%%%%%%%%%%%%%%%%%%%%%%%%%%%%%%%%%%%%%%%%
\begin{cor}\label{Xuniq}
Assume Hypothesis  \ref{Hmatr}. Suppose that $\Sigma$ and
$\widehat{\Sigma}$ are disjoint Borel subsets of $\bbR$
such that $\dist(\Sigma,\widehat{\Sigma})>0$. Let
$\cX=\cX(A,B,\Sigma,\widehat{\Sigma})$ be the set of all
bounded operators $X$ from $\cH_A$ to $\cH_C$ with the
properties
\begin{align}
\label{sZA}
&\spec(A+BX)\subset\Sigma ,\\
\label{sZC} &\spec(C-B^*X^*)\subset\widehat{\Sigma}, \quad
\dom(C-B^*X^*)=\dom(C).
\end{align}
Then if $X,Y\in\cX$ satisfy the Riccati equation
\eqref{Ric0}, then $X=Y$.
\end{cor}
%%%%%%%%%%%%%%%%%%%%%%%%%%%%%%%%%%%%%%%%%%%%%%%%%%%%%%%%%%
\begin{proof}
Suppose $X$ and $Y$ are two bounded solutions to
\eqref{Ric0} both satisfying \eqref{sZA} and \eqref{sZC}.
Then by Theorem \ref{thHi2} the graphs of $X$ and $Y$ both
coincide with the spectral subspace of the $2\times2$
operator matrix \eqref{twochannel} associated with the set
$\Sigma$, and hence, $X=Y$.
\end{proof}
%%%%%%%%%%%%%%%%%%%%%%%%%%%%%%%%%%%%%%%%%%%%%%%%%%%%%%%%

Under Hypothesis \ref{twochannel} introduce the
operator-valued Herglotz function
\begin{equation}
\label{M1def} M(\lambda)=\lambda I -A+B(C-\lambda
I)^{-1}B^*, \quad \lambda\in\rho(C).
\end{equation}

By definition the spectrum $\spec(M)$ of the function $M$
 is the set of all $\lambda\in\bbC$ such that either
the operator $M(\lambda)$ is not invertible or the inverse
$[M(\lambda)]^{-1}$ is an unbounded operator.

It is well known (see, e.g., \cite{MeMo99}) that the
resolvent of the operator $H$ \eqref{twochannel} can be
represented as the following $2\times2$ operator matrix
\begin{align}
\nonumber (H-\lambda I)^{-1}=&\begin{pmatrix} 0 &   0 \\ 0
& (C-\lambda I)^{-1}
\end{pmatrix}\\
\label{ResH}
&-\left(\begin{array}{c} I \\
-(C-\lambda I)^{-1}B^*
\end{array}\right)
M(\lambda)^{-1}
\begin{pmatrix}
I &\,\,  -B(C-\lambda I)^{-1}
\end{pmatrix}, \\
\nonumber & \quad  \lambda\in\rho(H),
\end{align}
where $M$ is the Herglotz function given by \eqref{M1def}.
Representation \eqref{ResH} shows that for
$\lambda\in\rho(C)$ the operator $H-\lambda I$ has a
bounded inverse iff $M(\lambda)$ does, which means that
\begin{equation}
\label{HM0} \spec(H)\cap\rho(C)=\spec(M)\cap\rho(C).
\end{equation}

We will also need the following general result (cf.
\cite[Theorem 2.2]{AdLT}, \cite[Proposition 2.4 and
Theorem 2.5]{MenShk}, and \cite[Theorems 4.4 and
5.1]{LMMT}).

%%%%%%%%%%%%%%%%%%%%%%%%%%%%%%%%%%%%%%%%%%%%%%%%%%%%%%%%%%%%%%%%%%%%%%%
\begin{thm}
\label{RicSol} Assume Hypothesis \ref{Hmatr} and suppose
that the Herglotz function \eqref{M1def} admits the
factorization
\begin{equation}
\label{M1fact} M(\lambda)=W(\lambda)(Z-\lambda
I),\quad\lambda\in\Omega,
\end{equation}
where $Z$ is a bounded operator in $\cH_A$ such that
$\spec(Z)\cap\spec(C)=\emptyset$, $\Omega$ is a domain in
$\rho(C)$ such that $\spec(Z)\subset\Omega$, and $W$ is a
holomorphic $\cB(\cH_A)$-valued function on $\Omega$, such
that for any $\lambda\in\spec(Z)$ the operator
$W(\lambda)$ has a bounded inverse. Then $\spec(Z)$ is an
isolated part of the spectrum of the operator $H$ and the
spectral subspace $\Ran\EE_H(\sigma(Z))$ of $H$ associated
with the set $\spec(Z)$ is the graph of the  bounded
operator $X$ from $\cH_A$ to $\cH_C$ given by
\begin{equation}
\label{Xdeff} X=-\frac{1}{2\pi{\mathrm i}} \int_\Gamma
d\lambda(C-\lambda I)^{-1}B^* (Z-\lambda I)^{-1}.
\end{equation}
Here, $\Gamma$ is an arbitrary Jordan contour in
$\rho(Z)\cap\rho(C)$ (maybe consisting of several simple
Jordan contours) encircling $\spec(Z)$ in the clockwise
direction and having winding number $0$ with respect to
the spectrum of $C$.

Moreover, the operator $Z$ can be written in terms of the
operator $X$ given by \eqref{Xdeff} as
\begin{equation}
\label{ZvX} Z=A+BX
\end{equation}
and the factor $W(\lambda)$ admits analytic continuation
to the whole resolvent set of the operator $C$ by the
following formula
\begin{equation}
\label{W} W(\lambda)=I-B(C-\lambda)^{-1}X, \quad
\lambda\in\rho(C).
\end{equation}
\end{thm}
%%%%%%%%%%%%%%%%%%%%%%%%%%%%%%%%%%%%%%%%%%%%%%%%%%%%%%%%%%%%%%%%%%%%%%%
%%%%%%%%%%%%%%%%%%%%%%%%%%%%%%%%%%%%%%%%%%%%%%%%%%%%%
\begin{proof}
By hypothesis the function $W(\lambda)$ is holomorphic on
an open set $\Omega\subset\rho(C)$ containing the closed
subset $\spec(Z)$ and the operator $W(\lambda)$ has a
bounded inverse for any $\lambda\in\spec(Z)$. Hence there
is an open neighborhood $\wOm$ of $\spec(Z)$ in $\Omega$
where the operator $W(\lambda)$ is boundedly invertible,
i.e. $W(\lambda)^{-1}\in\cB(\cH_A)$ for any
$\lambda\in\wOm$. By \eqref{M1fact}
\begin{equation}
\label{MZinv} M(\lambda)^{-1}=(Z-\lambda
I)^{-1}W(\lambda)^{-1}, \quad \lambda\in
\wOm\setminus\spec(Z).
\end{equation}
Taking into account \eqref{HM0} one infers that the
spectrum of $H$ in $\wOm$ coincides with that of $M$ and,
thus, with that of $Z$, that is,
\begin{equation}
\label{HUZ} \spec(H)\cap\wOm=\spec(Z),
\end{equation}
Since $\wOm$ is an open set, $\spec(Z)$ is a closed set,
and $\spec(Z)\subset\wOm$, one also concludes that
$\spec(Z)$ is isolated from the remaining part of the
spectrum of $H$, i.e.
\begin{equation*}
\dist(\spec(Z),\spec(H)\setminus\spec(Z))>0.
\end{equation*}

Using representation \eqref{ResH}, for the spectral
projection $\sE_H\bigl(\spec(Z)\bigr)$ of the operator $H$
associated with the set $\spec(Z)$ the Riesz integration
yields
\begin{align*}
\sE_H\bigl(\spec(Z) \bigr)=&
\sE_H\bigl(\spec(Z)\bigr)=\frac{1}{2\pi\ri}\int\limits_{\Gamma}
d\lambda(H-\lambda I)^{-1}\\
 =&-\frac{1}{2\pi\ri}\int\limits_{\Gamma}
d\lambda\begin{pmatrix} I \\
-(C-\lambda I)^{-1}B^*
\end{pmatrix}
M(\lambda)^{-1}
\begin{pmatrix}
I & \,\quad -B(C-\lambda I)^{-1}
\end{pmatrix},
\end{align*}
where $\Gamma$ stands for an arbitrary Jordan contour
(possibly consisting of several simple Jordan contours) in
$\wOm$ encircling the spectrum of $Z$ in the clockwise
direction and having winding number $0$ with respect to
the spectrum of $C$. Hence,
\begin{equation}
\label{EH} \sE_H\bigl(\spec(Z)\bigr)=
\begin{pmatrix}
E & G^*\\
G & F\\
\end{pmatrix},
\end{equation}
where
\begin{align}
\label{Eint}
E&=-\frac{1}{2\pi\ri}\int_\Gamma d\lambda M(\lambda)^{-1}, \\
\nonumber F&=-\frac{1}{2\pi\ri}\int_\Gamma
d\lambda(C-\lambda I)^{-1}B^* M(\lambda)^{-1}B(C-\lambda
I)^{-1},
\end{align}
and
\begin{align*}
\nonumber G&=\frac{1}{2\pi\ri}\int_\Gamma
d\lambda(C-\lambda I)^{-1}B^* M(\lambda)^{-1}\\
&=\frac{1}{2\pi\ri}\int_\Gamma d\lambda (C-\lambda
I)^{-1}B^* (Z-\lambda I)^{-1} W(\lambda)^{-1},
\end{align*}
using factorization formula \eqref{M1fact}. Since both the
operator-valued functions $(C-\lambda)^{-1}$ and
$W(\lambda)^{-1}$ are holomorphic in $\wOm$, applying the
Daletsky-Krein formula  (see \cite[Lemma I.2.1]{DK}) we
get
\begin{align}
\nonumber G=&\left[\frac{1}{2\pi\ri} \int_\Gamma d\lambda
(C-\lambda I)^{-1}B^*(Z-\lambda I)^{-1}\right]\\
\label{G} &\times \left[\frac{1}{2\pi\ri}\int_\Gamma
d\lambda (Z-\lambda I)^{-1}W(\lambda)^{-1}\right].
\end{align}
Hence, combining \eqref{MZinv} and \eqref{Eint} proves the
representation
\begin{equation}
\label{XE} G=XE,
\end{equation}
where $X$ (the first factor on the r.h.s. part of
\eqref{G}) is given by
\begin{equation}
\label{XinT} X=-\frac{1}{2\pi{\mathrm i}} \int_\Gamma
d\lambda(C-\lambda I)^{-1}B^* (Z-\lambda I)^{-1}.
\end{equation}
In an analogous way one also proves that
\begin{equation}
\label{XEX} F=XEX^*.
\end{equation}

Clearly, for $\lambda\in\Gamma$ we have $\Ran(C-\lambda
I)^{-1}\subset\dom(C)$ and hence
\begin{equation*}
\Ran(X)\subset\Dom(C),
\end{equation*}
which immediately follows from  \eqref{XinT}. Multiplying
both sides of \eqref{XinT} by $B$ from the left yields
\begin{equation} \label{BX1}
BX=-\frac{1}{2\pi{\mathrm i}}\int_\Gamma d\lambda
B(C-\lambda I)^{-1}B^*(Z-\lambda I)^{-1}.
\end{equation}
Meanwhile,
\begin{equation*}
B(C-\lambda I)^{-1}B^*=A-\lambda I+M(\lambda)= A-\lambda I
+W(\lambda)(Z-\lambda I), \quad \lambda\in\Gamma,
\end{equation*}
and, hence, using \eqref{BX1}
\begin{align}
\label{BXint} BX&=-\frac{1}{2\pi{\mathrm i}}\int_\Gamma
d\lambda
\left[W(\lambda)+A(Z-\lambda I)^{-1}-\lambda(Z-\lambda
I)^{-1}\right].
\end{align}
The function $W(\lambda)$ is holomorphic in the domain
bounded by the contour $\Gamma$ and, thus, the first term
in the integrand on the r.h.s.\ of \eqref{BXint} gives no
contribution. Since $\Gamma$ encircles the spectrum of
$Z$, the integration of the remaining two terms in
\eqref{BXint} can be performed explicitly using the
operator version of the residue theorem, which yields
$BX=-A+Z$ and hence
\begin{equation}
\label{Xgr} Z=A+BX
\end{equation}
proving  \eqref{ZvX}.

Since the spectra of the operators $C$ and $Z$ are
disjoint and $Z$ is a bounded operator, it is
straightforward to show (see, e.g., \cite{D53} or
\cite{R56})  that the operator $X$ given by \eqref{XinT}
is the unique solution to the operator Sylvester equation
\begin{equation*}
XZ-CX=B^*,
\end{equation*}
which by \eqref{Xgr} proves that $X$ solves the Riccati
equation \eqref{Ric0}.

Now applying Theorem \ref{thHi2} and using \eqref{EH},
\eqref{XE}, and \eqref{XEX} we arrive at the series of
equalities
\begin{align}
\nonumber I&= \frac{1}{2\pi\ri}\int_\Gamma
d\lambda(Z-\lambda I)^{-1}\\
\nonumber &=\frac{1}{2\pi\ri}\int_\Gamma d\lambda{\sP}_{A}
\begin{pmatrix}
(Z-\lambda I)^{-1} & 0 \\ 0 & (\widehat{Z}-\lambda I)^{-1}
\end{pmatrix}
{\sP}_{A}^*\\
\nonumber &=\frac{1}{2\pi\ri}\int_\Gamma d\lambda
{\sP}_{A} V^{-1}(H-\lambda I)^{-1}V
{\sP}_{A}^*\\
\nonumber &={\sP}_{A} V^{-1} \sE_H\bigl(\spec(Z)\bigr)
V{\sP}_{A}^*\\
\label{IXX} &=E(I+X^*X),
\end{align}
where $\widehat{Z}=C-B^*X^*$ with
$\Dom(\widehat{Z})=\Dom(C)$,
\begin{equation*}
V=\begin{pmatrix}
I & -X^* \\
X &  I
\end{pmatrix},
\end{equation*}
and ${\sP}_A$ is the canonical projection from $\cH$ to
$\cH_A$ (i.e., $\sP_A(f_A\oplus f_C)=f_A$ for
$f_A\in\cH_A$ and $f_C\in\cH_C$). Combining \eqref{XE},
and  \eqref{IXX} one concludes that the spectral
projection $\sE_{H}(\spec(Z))$ admits the representation
\begin{equation}
\label{sEQ} \sE_{H}(\spec(Z))=
\begin{pmatrix}
(I+X^*X)^{-1}  & (I+X^*X)^{-1}X^* \\
X(I+X^*X)^{-1} & X(I+X^*X)^{-1}X^*
\end{pmatrix}.
\end{equation}

Note that the contour $\Gamma$ in \eqref{XinT} can be
replaced by an arbitrary Jordan contour in
$\rho(Z)\cap\rho(C)$ (possibly consisting of several
simple Jordan contours) encircling the set $\spec(Z)$ in
the clockwise direction and having winding number $0$ with
respect to $\spec(C)$. Then observing that the r.h.s.\ of
\eqref{sEQ} is nothing but  the orthogonal projection in
$\cH=\cH_A\oplus\cH_C$ onto the graph of the operator $X$
proves that $\Ran\EE_H(\spec(Z))=\cG(X)$.

We conclude with the proof of the representation
\eqref{W}. First, by \eqref{Xgr} we notice that
\begin{align*}
[I-B(C-\lambda)^{-1}X](Z-\lambda)&=[I-B(C-\lambda)^{-1}X](A+BX-\lambda)\\
&=A+BX-\lambda -B(C-\lambda)^{-1}(XA+XBX-\lambda X), \\
&\quad \lambda\in\rho(C).
\end{align*}
Since $X$ solves the Riccati equation \eqref{Ric0} one
infers that $XA+XBX=B^*+CX$ which implies
\begin{align}
\nonumber
[I-B(C-\lambda)^{-1}X](Z-\lambda)&=M(\lambda)+BX-B(C-\lambda)^{-1}(C-\lambda)X\\
\label{dopM}
 &=M(\lambda), \quad \lambda\in\rho(C).
\end{align}
Hence combining \eqref{M1fact} and \eqref{dopM} yields
\begin{equation*}
W(\lambda)=M(\lambda)(Z-\lambda)^{-1}=I-B(C-\lambda)^{-1}X,
\quad \lambda\in\Omega\setminus\spec(Z).
\end{equation*}
Then the analytic continuation completes the proof.
\end{proof}
%%%%%%%%%%%%%%%%%%%%%%%%%%%%%%%%%%%%%%%%%%%%%%%%%%%%%%%%%%%%%%%%%%%%%%%
\begin{rem}
\label{RemAnal} By analytic continuation of both parts of
\eqref{M1fact} one concludes that the factorization
formula \eqref{M1fact} holds for any $\lambda\in\rho(C)$
with $W(\lambda)$ given by \eqref{W}. Moreover, the
representation
\begin{equation*}
[W(\lambda)]^{-1}=(Z-\lambda I)[M(\lambda)]^{-1}, \quad
\lambda\in\wOm\setminus\spec(M)=\wOm\setminus\spec(Z),
\end{equation*}
and the fact that $\spec(M)\subset\spec(H)$ imply that
$[W(\lambda)]^{-1}$ admits analytic continuation as a
$\cB(\cH_A)$-valued holomorphic function to the domain
$\rho(H)\cup\spec(Z)$.
\end{rem}

%%%%%%%%%%%%%%%%%%%%%%%%%%%%%%%%%%%%%%%%%%%%%%%%%%%%%%%%%%%%%%%%%%%%%%

%%%%%%%%%%%%%%%%%%%%%%%%%%%%%%%%%%%%%%%%%%%%%%%%%%%%%%%%%%%%%%%%%%%%%%%%%%
\section{An existence result. Proof of Theorem 1 (i)}
\label{SecExist}
%%%%%%%%%%%%%%%%%%%%%%%%%%%%%%%%%%%%%%%%%%%%%%%%%%%%%%%%%%%%%%%%%%%%%%%%%%
As we have already mentioned in Introduction  our main
technical tool in proving the solvability of the Riccati
equation
 is the Virozub-Matsaev factorization theorem \cite{ViMt}
(also see \cite{MrMt}). For convenience of the reader we
reproduce the corresponding statement following
Propositions 1.1 and 1.2 in \cite{MenShk}.

%%%%%%%%%%%%%%%%%%%%%%%%%%%%%%%%%%%%%%%%%%%%%%%%%%%%%%%%%%%%%%%%%%%%%%%
\begin{thm}
\label{Matsa} Let $\cK$ be a Hilbert space and
$F(\lambda)$ a holomorphic $\cB(\cK)$-valued function on a
simply connected domain $\Omega\subset\bbC$. Assume that
$\Omega$ includes an interval $[a,b]\subset\bbR$ such that
\begin{equation*}
F(a)< 0, \quad F(b)>
0,\quad\text{and}\quad\frac{d}{d\lambda}F(\lambda)> 0\quad
\text{for all}\quad\lambda\in[a,b].
\end{equation*}
Then there exist  a domain
$\widetilde{\Omega}\subset\Omega$ containing $[a,b]$ and a
unique bounded operator $Z$ on $\cK$ with
$\spec(Z)\subset(a,b)$ such that $F(\lambda)$ admits the
factorization
\begin{equation*}
F(\lambda)=G(\lambda)(Z-\lambda I), \quad
\lambda\in\widetilde{\Omega},
\end{equation*}
where $G(\lambda)$ is a holomorphic operator-valued
function on $\widetilde{\Omega}$ whose values are bounded
and boundedly invertible operators in $\cK$, that is,
\begin{equation*}
G(\lambda)\in\cB(\cK)\quad\text{and}\quad
[G(\lambda)]^{-1}\in\cB(\cK)\quad
\quad\lambda\in\widetilde{\Omega}.
\end{equation*}
\end{thm}
%%%%%%%%%%%%%%%%%%%%%%%%%%%%%%%%%%%%%%%%%%%%%%%%%%%%%%%%%%%%%%%

Based on Theorem \ref{Matsa} we obtain the following
factorization result, the  cornerstone for our further
considerations.

%%%%%%%%%%%%%%%%%%%%%%%%%%%%%%%%%%%%%%%%%%%%%%%%%%%%%%%%%%%
\begin{thm}
\label{Zexist} Assume Hypothesis \ref{Hmatr}. Assume, in
addition, that $C$ has a finite spectral gap $\Delta
=(\alpha,\beta)$, $\alpha<\beta$, the spectrum of $A$ lies
in $\Delta$, i.e., $\spec(A)\subset \Delta$, and
\begin{equation*}
\|B\|<\sqrt{d|\Delta|},
\end{equation*}
where
\begin{equation*}
d=\dist\bigl(\spec(A),\spec(C)\bigr).
\end{equation*}
Then there is a unique operator $Z\in\cB(\cH_A)$ with
$\spec(Z)\subset\Delta$ such that the operator-valued
function $M(\lambda)$ given by \eqref{M1def} admits the
factorization
 \begin{equation}\label{fact}
M(\lambda)=W(\lambda)(Z-\lambda I),\quad\lambda\in\rho(C),
\end{equation}
with a holomorphic $\cB(\cH_A)$-valued function $W$ on
$\rho(C)$. Moreover, for any
\begin{equation*}
\lambda\in\bigl(\bbC\setminus \spec(M)\bigr)\cup\Delta
\end{equation*}
the operator $W(\lambda)$ has a bounded inverse and
\begin{equation}\label{inclusion}
\spec(Z)=\spec(H)\cap\Delta
\subset[\inf\spec(A)-\delta_-,\sup\spec(A)+\delta_+],
\end{equation}
where
\begin{align}
\label{dBm}
\delta_-&=\|B\|\tan\left(\frac{1}{2}\arctan\frac{2\|B\|}
{\beta-\inf\spec(A)}\right)<\inf\spec(A)-\alpha ,\\
\label{dBp}
\delta_+&=\|B\|\tan\left(\frac{1}{2}\arctan\frac{2\|B\|}
{\sup\spec(A)-\alpha }\right)<\beta-\sup\spec(A).
\end{align}
\end{thm}
%%%%%%%%%%%%%%%%%%%%%%%%%%%%%%%%%%%%%%%%%%%%%%%%%%%%%%%%%%%
\begin{proof}
By the spectral theorem
\begin{equation*}
B(C-\lambda I)^{-1}B^*=\int\limits_{\R\setminus \Delta}
B\EE_C(d\mu)B^*\frac{1}{(\mu-\lambda)},\quad\lambda\in\rho(C),
\end{equation*}
where $\EE_C(\mu)$ stands for the spectral family of the
self-adjoint operator $C$. Hence
\begin{equation}
\label{DerE} \displaystyle\frac{d}{d\lambda}M(\lambda)=I +
\int\limits_{\R\setminus \Delta}
B\EE_C(d\mu)B^*\frac{1}{(\mu-\lambda)^2},
\qquad\lambda\in\rho(C).
\end{equation}
For $\lambda\in\Delta$ the integral in \eqref{DerE} is a
non-negative operator. Therefore, the derivative of
$M(\lambda)$ is a strictly positive operator:
\begin{equation*}
\label{Der}
 \displaystyle\frac{d}{d\lambda}M(\lambda)\geq I>0,
 \quad \lambda\in \Delta.
\end{equation*}

Next we estimate the quadratic form of $M(\lambda)$. Let
$f\in\cH_A$,  $\|f\|=1$. Then
\begin{align}\label{Fflam}
\lal M(\lambda)f,f\ral&=
   \lambda-\lal A f,f\ral+\lal(C-\lambda)^{-1}B^*f,B^*f\ral
\nonumber\\
&=\lambda-\lal A f,f\ral
 +\int\limits_{-\infty}^{\alpha }
\frac{1}{\mu-\lambda}\lal \EE_C(d\mu)B^*f,B^*f\ral \\
&\quad+\int\limits_{\beta}^{+\infty}
\frac{1}{\mu-\lambda}\lal \EE_C(d\mu)B^*f,B^*f\ral,
\quad\lambda\in\rho(C). \nonumber
\end{align}
Since for $\lambda\in \Delta$ the integral in the second
line of \eqref{Fflam} is non-positive and the one in the
third line  is non-negative, one obtains the two-sided
estimate
\begin{equation*}
\bigg ( \lambda-\sup
\spec(A)-\frac{\|B\|^2}{\lambda-\alpha }\bigg )I \le
 M(\lambda)
\le \bigg ( \lambda-\inf
\spec(A)-\frac{\|B\|^2}{\lambda-\beta} \bigg )I,
\quad\lambda\in \Delta.
\end{equation*}
 Now, a simple calculation  shows that
\begin{equation}
\label{Ier} M(\lambda)< 0 \quad \mbox{for}\quad
\lambda\in\bigl(\alpha ,\inf\spec(A)-\delta_-\bigr)
\end{equation}
and
\begin{equation}
\label{Iel} M(\lambda)> 0 \quad \mbox{for}\quad
\lambda\in\bigl(\sup\spec(A)+\delta_+,\beta\bigr),
\end{equation}
where $\delta_-$ and $\delta_+$ are given by \eqref{dBm}
and \eqref{dBp}, respectively.

Thus, the function $F(\lambda)=M(\lambda)$ satisfies
assumptions of Theorem \ref{Matsa} for any $a\in(\alpha
,\inf\spec(A)-\delta_-)$ and any
$b\in(\sup\spec(A)+\delta_+, \beta)$, proving the
existence of the unique bounded operator $Z\in\cB(\cH_A)$
such that \eqref{fact} and \eqref{inclusion} hold taking
into account \eqref{HM0}. It follows from Theorem
\ref{Matsa} that the factor $W(\lambda)$ in \eqref{M1fact}
has a bounded inverse in a complex neighborhood
$U\subset\bbC$ of the interval
$[\inf\spec(A)-\delta_-,\sup\spec(A)+\delta_+]$. Moreover,
the operator $W(\lambda)$ has a bounded inverse for any
$\lambda\in\bigl(\bbC\setminus \spec(M)\bigr)\cup\Delta$
by Remark \ref{RemAnal}. The proof is complete.
\end{proof}
%%%%%%%%%%%%%%%%%%%%%%%%%%%%%%%%%%%%%%%%%%%%%%%%%%%%%%%%%%%%%%%%%%%%

We are ready to prove the part (i) of Theorem 1.

\begin{proof}[Proof of Theorem \ref{thm:1} (i)]
By Theorem \ref{Zexist} the Herglotz function
\eqref{M1def} admits the factorization \eqref{M1fact} with
$W(\lambda)$ and $Z$ satisfying hypothesis of Theorem
\ref{RicSol}. Therefore, the Riccati equation \eqref{Ric0}
has a bounded solution $X$ given by \eqref{Xdeff}. Theorem
\ref{RicSol} also shows that the graph $\cG(X)$ of the
operator $X$ coincides with the spectral subspace $\Ran
\sE_H\bigl(\spec(Z)\bigr)$ for the block operator matrix
$H$ given by \eqref{twochannel}. By \eqref{inclusion} the
subspace $\Ran \sE_H\bigl(\spec(Z)\bigr)$ coincides with
$\Ran \sE_H\bigl(\Delta\bigr)$, proving that
\begin{equation}
\label{EHG} \Ran\sE_H\bigl(\Delta\bigr)=\cG(X).
\end{equation}

To prove that $X$ possesses the properties \eqref{Uniq} we
proceed as follows. Let
\begin{equation}
V=\left(\begin{array}{cc}
I & -X^* \\
X & I
\end{array}\right), \quad S=(I+X^*X)^{1/2}, \quad\text{and}\quad
\widehat{S}=(I+XX^*)^{1/2}.
\end{equation}
{}From Theorem \ref{thHi2} it follows that
\begin{equation}
\label{EHZZ} \EE_H(\Delta)=V\left(\begin{array}{cc}
S^{-1}\EE_\Lambda(\Delta)S & 0 \\
0 &
\widehat{S}^{-1}\EE_{\widehat{\Lambda}}(\Delta)\widehat{S}
\end{array}\right)V^{-1},
\end{equation}
where $\Lambda$ and $\widehat{\Lambda}$ are the
self-adjoint operators defined by \eqref{HA} and
\eqref{HC}, respectively; $\EE_\Lambda(\Delta)$ and
$\EE_{\widehat{\Lambda}}(\Delta)$ denote the spectral
projections for $\Lambda$ and $\widehat{\Lambda}$
associated with the interval $\Delta$. Since $\Lambda$ is
similar to $Z$ ($\Lambda=SZS^{-1}$) and by Theorem
\ref{Zexist} the inclusion $\spec(Z)\subset\Delta$ holds,
one concludes that $\spec(\Lambda)\subset\Delta$. Hence,
$\EE_\Lambda(\Delta)=I$ and then \eqref{EHZZ} implies the
equality
\begin{equation}
\label{EHZZ1} \EE_H(\Delta)= V \left( \begin{array}{cc}
I & 0 \\
0 & 0
\end{array}\right) V^{-1}
+ V\left(\begin{array}{cc}
0 & 0 \\
0 &
\widehat{S}^{-1}\EE_{\widehat{\Lambda}}(\Delta)\widehat{S}
\end{array}\right)V^{-1}.
\end{equation}
The first summand on the r.h.s.\ of \eqref{EHZZ1}
\begin{equation*}
V \left( \begin{array}{cc}
I & 0 \\
0 & 0
\end{array}\right) V^{-1}=\begin{pmatrix}
(I+X^*X)^{-1}  & (I+X^*X)^{-1}X^* \\
X(I+X^*X)^{-1} & X(I+X^*X)^{-1}X^*
\end{pmatrix}
\end{equation*}
coincides with the orthogonal projection onto the graph
subspace $\cG(X)$ and, hence, by \eqref{EHG} it equals
$\sE_H\bigl(\Delta\bigr)$. Thus, the second summand on the
r.h.s.\ of \eqref{EHZZ1} vanishes,
\begin{equation*}
V\left(\begin{array}{cc}
0 & 0 \\
0 &
\widehat{S}^{-1}\EE_{\widehat{\Lambda}}(\Delta)\widehat{S}
\end{array}\right)V^{-1}=0,
\end{equation*}
which means that $\EE_{\widehat{\Lambda}}(\Delta)=0$ and,
therefore, $\spec(\widehat{\Lambda})\cap\Delta=\emptyset$.
By \eqref{HC} it is now straightforward to see that $X$
satisfies \eqref{Uniq}. It remains to conclude that by
Corollary \ref{Xuniq} the operator $X$ is the unique
bounded solution to \eqref{Ric0} satisfying \eqref{Uniq}
which completes the proof of Theorem 1 (i).
\end{proof}

%%%%%%%%%%%%%%%%%%%%%%%%%%%%%%%%%%%%%%%%%%%%%%%%%%%%%%%%%%%%%
\begin{rem}
\label{coptim} Obviously, under hypothesis of Theorem
\ref{thm:1} the inequality $|\Delta|\geq 2d$ holds, with
the equality sign occurring only if the spectrum of the
operator $A$ is a one point set. Hence,  the  condition
$\|B\|<\sqrt{2}d$, which is stronger than \eqref{korint},
implies  the existence of a bounded  solution to the
Riccati equation \eqref{Ric0}. This means that in the case
where $C$ has a finite spectral gap $\Delta$ and
$\spec(A)\subset\Delta$ the best possible constant
$c_{\mathrm{best}}$ in condition \eqref{cdest} ensuring
the solvability of the Riccati equation satisfies the
inequality $c_{\mathrm{best}}\geq \sqrt{2}$. On the other
hand in \cite[Lemma 3.11 and Remark 3.12]{AMM} it is shown
that $c_{\mathrm{best}}\leq \sqrt{2}$. Thus, $c=\sqrt{2}$
is best possible in inequality \eqref{cdest} which ensures
the solvability of the Riccati equation under the
additional hypothesis that $\dist(\sigma(A), \sigma
(C))>0$ and that the spectrum of $ A$ lies in a spectral
gap of the operator $C$.
\end{rem}

%%%%%%%%%%%%%%%%%%%%%%%%%%%%%%%%%%%%%%%%%%%%%%%%%%%%%%%%%%%%%%
\section{The $\text{tan}{\,\Theta}$ Theorem} \label{SecTan}
%%%%%%%%%%%%%%%%%%%%%%%%%%%%%%%%%%%%%%%%%%%%%%%%%%%%%%%%%%%%%%

We start out by recalling a concept of the operator angle
between two subspaces in a Hilbert space going back to the
works by Friedrichs \cite{Friedrichs}, M.~Krein,
Kransnoselsky, and Milman \cite{Krein:Krasnoselsky},
\cite{Krein:Krasnoselsky:Milman}, Halmos \cite{Halmos:69},
and Davis and Kahan \cite{Davis:Kahan}. A comprehensive
discussion of this concept can be found, e.g., in
\cite{KMMgeom}.

Given a closed subspace $\cQ$ of the Hilbert space
$\cH=\cH_A\oplus\cH_C$, introduce the operator angle
$\Theta$ between the subspaces $\cH_A\oplus\{0\}$ and
$\cQ$ by
\begin{equation}
\label{sinT}
\Theta=\arcsin\sqrt{I_{\cH_A}-{\sP}_A\sQ{\sP}_A^*},
\end{equation}
where ${\sP}_A$ is the canonical projection from $\cH$
onto $\cH_A$ and $\sQ$ the orthogonal projection in $\cH$
onto $\cQ$. If the subspace $\cQ$ is the graph $\cG(X)$ of
a bounded operator $X$ from $\cH_A$ to $\cH_C$, then (see
\cite{KMMgeom}; cf.\ \cite{Davis:Kahan} and
\cite{Halmos:69})
\begin{equation}
\label{tan-X} \tan\Theta=\sqrt{X^*X}
\end{equation}
and
\begin{equation}
\label{sinTPQ} \|\sin\Theta\|=\|\sQ-\sP\|,
\end{equation}
where $\sP={\sP_A}^{\!\!\!*}\sP_A$ denotes the orthogonal
projection in $\cH$ onto the subspace $\cH_A\oplus\{0\}$.

Note that the common definition of the operator angle
(see, e.g., \cite{KMMgeom}) slightly differs from
\eqref{sinT}. Usually, the operator angle is defined as
the restriction of the operator \eqref{sinT} onto the
maximal subspace of $\cH_A$ where it has a trivial kernel.
Clearly, the difference in these two definitions does not
effect the value of the norm $\|\tan\Theta\|$.

Now we are ready to prove the second principal result of
the paper, a generalization of the Davis-Kahan
$\tan\Theta$ Theorem \cite{Davis:Kahan}.

%%%%%%%%%%%%%%%%%%%%%%%%%%%%%%%%%%%%%%%%%%%%%%%%%%%%%%%%%%%%%%%%%%%%%%%%%%%%
\begin{proof}[Proof of Theorem \ref{thm:2}]
By hypothesis the Riccati equation \eqref{Ric0} has a
bounded solution $X$. Then, by Theorem \ref{thHi2}  the
operator $Z=A+BX$ is similar to the bounded self-adjoint
operator $\Lambda$ given by \eqref{HA} and hence
\begin{equation}
\label{spZA} \spec(Z)=\spec(\Lambda)\subset\bbR\,.
\end{equation}
The Riccati equation \eqref{Ric0} can be rewritten in the
form
\begin{equation}
\label{SylHelp} X(Z-\gamma I)-(C-\gamma I)X=B^*,
\end{equation}
where
\begin{equation}
\label{lamcen1}
\gamma=\frac{1}{2}\bigl(\sup\spec(Z)+\inf\spec(Z)\bigr)
=\frac{1}{2}\bigl(\sup\spec(\Lambda)+\inf\spec(\Lambda)\bigr).
\end{equation}
By hypothesis
$\dist\bigl(\spec(H|_{\cG(X)}),\spec(C)\bigr)=\delta>0$
and then Theorem \ref{thHi2} implies that
\begin{equation*}
\spec(H|_{\cG(X)})=\spec(Z)=\spec(\Lambda),
\end{equation*}
which by hypothesis proves the inclusion
\begin{equation}
\label{spC}
\spec(C)\subset\bigl(-\infty,\inf\spec(\Lambda)-\delta\bigr]\cup
\bigl[\sup\spec(\Lambda)+\delta,\infty\bigr).
\end{equation}
Hence, combining \eqref{spZA}, \eqref{lamcen1}, and
\eqref{spC} proves that $\gamma\in\rho(C)$ and
\begin{equation}
\label{CmA} \|(C-\gamma
I)^{-1}\|=\frac{1}{\|\Lambda-\gamma I\|+\delta}.
\end{equation}
Multiplying both sides of \eqref{SylHelp} by $(C-\gamma
I)^{-1}$ from the left one gets
 the
representation
\begin{equation}
\label{SylZC} X=(C-\gamma  I)^{-1}\bigl(X(Z-\gamma
I)-B^*\bigr).
\end{equation}
Using   Theorem \ref{thHi2} (ii) one obtains the estimate
\begin{align}
\nonumber \|X(Z-\gamma
I)\|=&\|X(I+X^*X)^{-1/2}(\Lambda-\gamma I)
(I+X^*X)^{1/2}\| \\
\label{NX1} &\leq \|X(I+X^*X)^{-1/2}\|\, \|\Lambda-\gamma
I\| (1+\|X\|^2)^{1/2}.
\end{align}
Clearly, by the spectral theorem
\begin{equation*}
\|X(I+X^*X)^{-1/2}\|=\sqrt{\|X^*X(I+X^*X)^{-1}\|}
=\frac{\|X\|}{(1+\|X\|^2)^{1/2}}.
\end{equation*}
Hence \eqref{NX1} implies the estimate
\begin{equation*}
\|X(Z-\gamma I)\|\leq \|X\|\, \|(\Lambda-\gamma I)\|,
\end{equation*}
which together with \eqref{SylZC} proves the norm
inequality
\begin{equation}
\label{XN2} \|X\|\leq \|(C-\gamma
I)^{-1}\|\bigl(\|\Lambda-\gamma I\|\|X\|+\|B\|\bigr).
\end{equation}
Solving inequality \eqref{XN2} with respect to $\|X\|$ and
taking into account \eqref{CmA} proves the  norm estimate
\eqref{NXest}. Finally, since
$\bigl\|\sqrt{X^*X}\bigr\|=\|X\|$, by the definition of
the operator angle \eqref{tan-X}
 one gets
\begin{equation*}
\|\tan\Theta\|=\|X\|.
\end{equation*}
Hence, \eqref{NXest} is equivalent to \eqref{TanTheta}.
\end{proof}
%%%%%%%%%%%%%%%%%%%%%%%%%%%%%%%%%%%%%%%%%%%%%%%%%%%%%%%%%%%%
\begin{rem}
It is natural to ask whether  estimate \eqref{TanTheta}
remains to hold if one replaces the distance
$\delta=\dist(\sigma(H|_{\cG(X)}),\sigma(C))$ by
\begin{equation*}
\widehat{\delta}=\dist\bigl(\spec(H|_{\cG(X)^\perp},\spec(A)\bigr)
=\dist\bigl(\spec(C-B^*X^*),\spec(A)\bigr).
\end{equation*}
The answer is ``No'': Example 6.1 in \cite{Davis:Kahan}
shows that the inequality
$\widehat{\delta}\|\tan\Theta\|\leq\|B\|$ fails to hold in
general.
\end{rem}

%%%%%%%%%%%%%%%%%%%%%%%%%%%%%%%%%%%%%%%%%%%%%%%%%%%%%%%%%%%%%%%%%%%%%%%%%%%%%
\section{Norm Estimates of Solutions. Proof~of~Theorem~1 (ii)}
\label{SecEst}
%%%%%%%%%%%%%%%%%%%%%%%%%%%%%%%%%%%%%%%%%%%%%%%%%%%%%%%%%%%%%%%%%%%%%%%%%%%%%

The  existence result of Theorem 1 (i) by itself gives no
clue for  estimating the norm of the corresponding
solution $X$ to the Riccati equation. To the contrary,
Theorem 2 provides such an estimate whenever some
additional information on the spectrum location of the
``perturbed'' operator matrix $H=\begin{pmatrix} A&
B\\B^*&C\end{pmatrix}$ is available. In turn, the  bounds
on the spectrum of the part $H|_{\cG(X)}$ of the operator
matrix $H$ associated with the reducing subspace $\cG(X)$
(needed to satisfy the hypotheses of Theorem \ref{thm:2})
can be obtained by combining the results of Theorems
\ref{RicSol} and \ref{Zexist}. As a result of performing
this program one gets an \emph{a priori} estimate on the
norm of  the solution $X$.

%%%%%%%%%%%%%%%%%%%%%%%%%%%%%%%%%%%%%%%%%%%%%%%%%%%%%%%%%%%%
\begin{thm}
\label{TanFin} Assume hypothesis of Theorem \ref{thm:1}
(i) with $\Delta=(\alpha,\beta)$, $\alpha<\beta$. Let  $X$
be the unique solution to the Riccati equation
\eqref{Ric0} referred to in Theorem \ref{thm:1}. Then
\begin{equation}
\label{XestFin} \|X\|\leq\frac{\|B\|}{\widetilde{\delta}},
\end{equation}
where
\begin{equation}
\label{dtilde}
\widetilde{\delta}=\min\{\inf\spec(A)-\alpha-\delta_-,\,\,
\beta-\sup\spec(A)-\delta_+\}>0
\end{equation}
with $\delta_\pm$ given by \eqref{dBm} and \eqref{dBp}.

\end{thm}
%%%%%%%%%%%%%%%%%%%%%%%%%%%%%%%%%%%%%%%%%%%%%%%%%%%%%%%%
\begin{proof} Under
hypothesis of Theorem \ref{thm:1} (i) with
$\Delta=(\alpha,\beta)$, $\alpha<\beta$ one can apply
Theorem \ref{Zexist} and  using the same strategy of proof
as that of Theorem \ref{thm:1} (i) one concludes that
\begin{equation*}
Z=A+BX,
\end{equation*}
 where
  $Z$ is the unique
operator with $\spec(Z)\subset\Delta$ referred to in
Theorem \ref{Zexist} and
 $X$ is the unique solution to the Riccati equation
\eqref{Ric0} referred to in Theorem \ref{thm:1} (i).
Hence, by Theorem \ref{Zexist}
\begin{equation}
\label{2st}
\spec(A+BX)\subset[\inf\spec(A)-\delta_-,\sup\spec(A)+\delta_+].
\end{equation}
By hypothesis (of Theorem \ref{thm:1} (i))
\begin{equation}
\label{3st} \spec(C)\subset\bigl(-\infty,\alpha \bigr]\cup
\bigl[\beta,\infty\bigr),
\end{equation}
which together with  \eqref{2st} yields the inequality
\begin{equation*}
\dist\bigl(\spec(A+BX),\spec(C)\bigr)\geq\widetilde{\delta}.
\end{equation*}
By Theorem \ref{RicSol} one observes that
$\spec(H)|_{\cG(X)}=\spec(Z)$ which proves \eqref{XestFin}
using Theorem \ref{thm:2}.
\end{proof}
%%%%%%%%%%%%%%%%%%%%%%%%%%%%%%%%%%%%%%%%%%%%%%%%%%%%%%%%%%%%%%%%%%%%%%%%

The proof of  Theorem \ref{thm:1} (ii) needs complementary
considerations. Our reasoning is  based on the celebrated
Davis-Kahan $\tan 2 \Theta$-Theorem \cite{Davis:Kahan}.
For convenience of the reader we reproduce  the
corresponding result (cf.  Theorem 2.4 (iii) and Remark
2.8 in \cite{KMMalpha} and Corollary 6.4 in
\cite{KMMgeom}).

\begin{thm}
\label{DaKa} Assume Hypothesis \ref{Hmatr}. Suppose that
the operator $C$ is bounded and
$\sup\spec(A)<\inf\spec(C)$. Then the open interval
$(\sup\spec(A),\inf\spec(C))$ is a spectral gap of the
operator $H$ and the spectral subspace
$\EE_H\bigl((-\infty,\sup\spec(A)]\bigr)$ is the graph of
a contractive operator $X$ from $\cH_A$ to $\cH_C$.
Moreover, the operator $X$ is the unique contractive
solution to the Riccati equation \eqref{Ric0} and its norm
satisfies the estimate
\begin{equation*}
%\label{DaKaE}
\|X\|\le \tan \left(\frac{1}{2}\arctan
\frac{2\|B\|}{d}\right) < 1,
\end{equation*}
where
\begin{equation*}
d=\dist(\spec(A),\spec(C)).
\end{equation*}
\end{thm}

Now we are prepared to prove Theorem 1 (ii).

\begin{proof}[Proof of Theorem \ref{thm:1} (ii)]
Let $X$ be the solution to the Riccati equation
\eqref{Ric0} referred to in Theorem \ref{thm:1} (i) and
thus the spectral subspace of the operator $H$ associated
with the interval $\Delta$ is the graph $\cG(X)$ of the
operator $X$. Then $\cG(X)$ is also the spectral subspace
of the operator $(H-\gamma I)^2$ associated with the
interval $[0,|\Delta|^2/4)$ where $\gamma$ is the center
of the interval $\Delta$, that is,
\begin{equation}\label{kon}
\cG(X)=\Ran \big (\EE_H(\Delta)\big )=\Ran \big
(\EE_{(H-\gamma I)^2} [0,{|\Delta |^2}/{4})\big ).
\end{equation} By inspection
one obtains that with respect to the decomposition
$\cH=\cH_A\oplus\cH_A$ the non-negative operator
$(H-\gamma I)^2$ reads
\begin{equation}\label{kvadrat}
(H-\gamma I)^2=\begin{pmatrix}
\widehat{A} & \widehat{B}\  \\
\widehat{B}^* & \widehat{C}
\end{pmatrix},
\end{equation}
where $\widehat{A}=(A-\gamma I)^2+BB^*$,
$\widehat{B}=AB+BC$,\, and\, $\widehat{C}=(C-\gamma
I)^2+B^*B$.

The hypothesis that the spectrum of $C$ lies in
$\bbR\setminus\Delta$ implies the operator inequality
\begin{equation*}
\widehat{C}\geq \frac{|\Delta|^2}{4}I+BB^*\geq
\frac{|\Delta|^2}{4}\,I.
\end{equation*}
The hypothesis $\spec(A)\subset\Delta$ yields
\begin{equation}\label{puta}
0\le\widehat{A}\leq
\left(\frac{|\Delta|}{2}-d\right)^2I+BB^*\leq
\left[\left(\frac{|\Delta|}{2}-d\right)^2+\|B\|^2\right]I,
\end{equation}
taking into account that $\dist(\spec(A),\spec(C))=d$.
Hence, under hypothesis \eqref{korints}  one concludes
that
\begin{equation}\label{dura}
\dist\bigl(\spec(\widehat{A}),\spec(\widehat{C})\bigr)\geq
d(|\Delta|-d)-\|B\|^2>0
\end{equation}
and that the spectra $\spec(\widehat{A})$ and
$\spec(\widehat{C})$ of the entries $\widehat{A}$ and
$\widehat{C}$ are subordinated, that is,
$\sup\spec(\widehat{A})<\inf\spec(\widehat{C})$. By
Theorem \ref{DaKa}  (cf. Theorem 2.1 in \cite{AL95}) one
infers that the interval
$\bigl(\sup\spec(\widehat{A}),\inf\spec(\widehat{C})\bigr)$
lies in the resolvent set of the operator $(H-\gamma
I)^2$. In particular, the following inclusion holds
\begin{equation}
\label{gap2}
\bigl((|\Delta|/2-d)^2+\|B\|^2,|\Delta|^2/4\bigr)\subset
\rho\bigl((H-\gamma I)^2\bigr).
\end{equation}
Therefore, the spectral subspaces of the operator
$(H-\gamma I)^2$ associated with the intervals
$[0,|\Delta|^2/4)$ and $[0,|(\Delta|/2-d)^2+\|B\|^2]$,
respectively, coincide, that is,
\begin{equation}\label{ruba}
\begin{split}
\cG(X) &=\Ran \bigl(\EE_{(H-\gamma I)^2}\bigl([0,|\Delta|^2/4)\bigr)\bigr)\\
& =\Ran \bigl(\EE_{(H-\gamma
I)^2}([0,(|\Delta|/2-d)^2+\|B\|^2])\bigr).
\end{split}
\end{equation}
{}From \eqref{kvadrat} one concludes that the operator
matrix $(H-\gamma I)^2$ is an off-diagonal perturbation of
the matrix $\diag\{\widehat{A}, \widehat{C}\}$ diagonal
with respect to the decomposition $\cH=\cH_A\oplus \cH_C$.
Applying again Theorem \ref{DaKa} proves that the spectral
subspace  of the operator $(H-\gamma I)^2$ associated with
the interval $[0,|(\Delta|/2-d)^2+\|B\|^2]$ is the graph
of a contraction $\widehat{X}$ satisfying the
norm-estimate \eqref{XestFin2}. By \eqref{kon} and
\eqref{ruba}  this subspace coincides with $\cG(X)$ and,
therefore, $\widehat{X}=X$ and hence \eqref{XestFin2}
holds for $X$, completing the proof.
\end{proof}

%%%%%%%%%%%%%%%%%%%%%%%%%%%%%%%%%%%%%%%%%%%%%%%%%%%%%%%%%%%%%

\begin{rem}
\label{controptim} Condition \eqref{korints} ensuring the
strict contractivity of the solution $X$ is sharp. This
can be seen as follows. Let $\cH_A=\bbC$, $\cH_C=\bbC^2$,
$A=0$,
\begin{equation*}
C=\left(\begin{array}{rr}
-d & 0 \\
 0 & d
\end{array}\right),\quad d>0,
\end{equation*}
and
\begin{equation*}
B=\left(\frac{b}{\sqrt{2}},\frac{b}{\sqrt{2}}\right),
\quad b\in\bbR.
\end{equation*}
By inspection one proves that the $2\times1$ matrix
\begin{equation}
\label{Xexem} X=\left(\begin{array}{r}
-\frac{b}{\sqrt{2}\,d} \\[0.5em]
\frac{b}{\sqrt{2}\,d}
\end{array}\right)
\end{equation}
solves the Riccati equation
\begin{equation*}
XA-CX+XBX=B^*.
\end{equation*}
Moreover,
\begin{equation*}
A+BX=0
\end{equation*}
and $X$ possesses the properties \eqref{Uniq}. Clearly,
$\|B\|=b$ and $ \|X\|=\frac{b}{d}=\frac{\|B\|}{d}. $

For $b\in [d, \sqrt{2}d)$ hypothesis \eqref{korint} is
satisfied with $\Delta=(-d, d)$, condition \eqref{korints}
fails to hold, and $\|X\|\ge 1$, that is, estimate
\eqref{korints} is sharp.
\end{rem}

%%%%%%%%%%%%%%%%%%%%%%%%%%%%%%%%%%%%%%%%%%%%%%%%%%%%%%%%%%%%%%%%%%%%%%%%

\bigskip

\subsection*{Acknowledgments}

V.~Kostrykin in grateful to V.~Enss, A.~Knauf, H.~Leschke,
and R.~Schrader for useful discussions. K.~A.~Makarov is
indebted to Graduiertenkolleg ``Hierarchie und Symmetrie
in mathematischen Modellen" for kind hospitality during
his stay at RWTH Aachen in summer 2002. A.~K.~Motovilov
acknowledges the kind hospitality and support by the
Department of Mathematics, University of Missouri,
Columbia, MO, USA. He was also supported in part by the
Russian Foundation for Basic Research.

\bigskip


\begin{thebibliography}{00}

\bibitem{AL95} V.~M.~Adamjan and H.~Langer, \textit{Spectral
properties of a class of rational operator valued
functions}, J. Operator Theory \textbf{33} (1995), 259 --
277.

\bibitem{AdLT} V.~Adamjan, H.~Langer, and C.~Tretter,
\textit{Existence and uniqueness of contractive solutions
of some Riccati equations}, J. Funct. Anal. \textbf{179}
(2001), 448 -- 473.

\bibitem{AMM}
S.\,Albeverio, K.\,A.\,Makarov, and A.\,K.\,Motovilov,
\textit{Graph subspaces and the spectral shift function},
Canad. J. Math. (to appear); arXiv:\,math.SP/0105142.


\bibitem{D53} Yu.~L.~Dalecki\u{\i},  On the asymptotic
solution of a vector differential equation, \textit{Dokl.
Akad. Nauk SSSR} \textbf{92} (1953), 881 -- 884 (Russian).

\bibitem{DK} Ju.~L.~Dalecki\u{\i} and M.~G.~Kre\u{\i}n,
\textit{Stability of Solutions of Differential Equations
in Banach Spaces}, Translations of Mathematical
Monographs, Vol.\,43, AMS, Providence, Rhode Island, 1974.

\bibitem{Daughtry} J.~Daughtry, \textit{Isolated solutions of
quadratic matrix equations}, Linear Algebra Appl.
\textbf{21} (1978), 89 -- 94.

\bibitem{Davis:Kahan} C.~Davis and W.~M.~Kahan, \textit{The
rotation of eigenvectors by a perturbation. III}, SIAM J.
Numer. Anal. \textbf{7} (1970), 1 -- 46.

\bibitem{Friedrichs} K.~Friedrichs, \textit{On certain
inequalities and characteristic value problems for
analytic functions and for functions of two variables},
Trans. Amer. Math. Soc. \textbf{41} (1937), 321 -- 364.

\bibitem{Halmos:69} P.~R.~Halmos, \textit{Two subspaces}, Trans.
Amer. Math. Soc. \textbf{144} (1969), 381--389.

\bibitem{KMMgeom} V.~Kostrykin, K.~A.~Makarov, and A.~K.~Motovilov,
\textit{Existence and uniqueness of solutions to the
operator Riccati eqution. A geometric approach}, in
Yu.~Karpeshina, G.~Stolz, R.~Weikard, Y.~Zeng (Eds.),
\textit{Advances in Differential Equations and
Mathematical Physics}, Contemporary Mathematics
\textbf{327}, Amer. Math. Soc., 2003 (to appear);
arXiv:\,math.SP/0207125.

\bibitem{KMMalpha} {V.\,Kostrykin, K.\,A.\,Makarov, and A.\,K.\,Motovilov,
\textit{A generalization of the $\tan\,2\Theta$
the\-o\-rem},} arXiv: \htmladdnormallink{math.SP/}{http://www.arXiv.org/abs/math.SP/0302020}
\htmladdnormallink{0302020}{http://www.arXiv.org/abs/math.SP/0302020}.

\bibitem{KMMgamma} V.~Kostrykin, K.~A.~Makarov, and A.~K.~Motovilov,
\textit{Perturbation of spectra and spectral subspaces},
preprint (2003).

\bibitem{Krein:Krasnoselsky} M.~G.~Krein and M.~A.~Krasnoselsky,
\textit{Fundamental theorems about extensions of Hermite
operators and some applications to the theory of
orthogonal polynomials and to the moment problem}, Uspekhi
Mat. Nauk \textbf{2} (1947), 60 -- 106 (Russian).

\bibitem{Krein:Krasnoselsky:Milman} M.~G.~Krein,
M.~A.~Krasnoselsky, and D.~P.~Milman, \textit{On defect
numbers of linear operators in Banach space and some
geometric problems}, Sbornik Trudov Instituta Matematiki
Akademii Nauk Ukrainskoy SSR, \textbf{11} (1948), 97 --
112 (Russian).

\bibitem{LMMT} H.~Langer, A.~Markus, V.~Matsaev, and C.~Tretter,
\textit{A new concept for block operator matrices: the
quadratic numerical range}, Linear Algebra Appl.
\textbf{330} (2001), 89 -- 112.

\bibitem{MrMt} A.~S.~Markus and V.~I.~Matsaev,
\textit{Spectral theory of holomorphic operator-functions
in Hilbert space}, Funct. Anal. Appl. \textbf{9} (1975),
73 -- 74.

\bibitem{MeMo99} R.~Mennicken and A.~K.~Motovilov,
\textit{Operator interpretation of resonances arising in
spectral problems for $2\times2$ operator matrices}, Math.
Nachr. \textbf{201} (1999), 117 -- 181;
arXiv:\,funct-an/9708001.

\bibitem{MenShk} R.~Mennicken and A.~A.~Shkalikov,
\textit{Spectral decomposition of symmetric operator
matrices}, Math. Nachr. \textbf{179} (1996), 259 -- 273.

\bibitem{Motovilov:SPb:91} A.~K.~Motovilov, \textit{Potentials
appearing after the removal of energy-dependence and
scattering by them}, in Proc. Intern. Workshop
``Mathematical Aspects of the Scattering Theory and
Applications'', St.~Petersburg State University,
St.~Petersburg, 1991, pp.~101 -- 108.

\bibitem{Motovilov:95} A.~K.~Motovilov, \textit{Removal of the
resolvent-like energy dependence from interactions and
invariant subspaces of a total Hamiltonian}, J. Math.
Phys. \textbf{36} (1995), 6647 -- 6664;
arXiv:\,funct-an/9606002.

\bibitem{R56} M.~Rosenblum,  \textit{On the operator
equation $BX-XA=Q$}, Duke Math. J. \textbf{23} (1956), 263
-- 269.

\bibitem{ViMt} A.~I.~Virozub and V.~I.~Matsaev, \textit{The
spectral properties of a certain class of self-adjoint
operator functions}, Funct. Anal. Appl. \textbf{8} (1974),
1 -- 9.

\end{thebibliography}
\end{document}